\documentclass{amsart}

\usepackage{amsmath}
\usepackage{amsfonts}
\usepackage{amssymb}

\newcommand{\half}{\tfrac{1}{2}}
\newcommand{\summ}{\mathop{{\sum}^{\star}}}
\newcommand{\sumh}{\mathop{{\sum}^{h}}}

\numberwithin{equation}{section}

\newtheorem{theorem}[subsection]{Theorem}
\newtheorem{proposition}[subsection]{Proposition}
\newtheorem{lemma}[subsection]{Lemma}

\begin{document}

\title{On the subconvexity problem for $GL(3)\times GL(2)$ $L$-functions}
\author{Rizwanur Khan}
\address{Mathematisches Institut, Georg-August Universit\"{a}t G\"{o}ttingen, Bunsenstra$\ss$e 3-5,
         D-37073 G\"{o}ttingen, Germany}
\email{rrkhan@uni-math.gwdg.de}
\thanks{The author was supported by a grant from the European Research Council (grant agreement number 258713)}

\begin{abstract}
Fix $g$ a self-dual Hecke-Maass form for $SL_3(\mathbb{Z})$. Let $f$ be a holomorphic newform of prime level $q$ and fixed weight. Conditional on a lower bound for a short sum of squares of Fourier coefficients of $f$, we prove a subconvexity bound in the $q$ aspect for $L(s, g\times f)$ at the central point.

\end{abstract}

\maketitle

\section{Introduction}

An outstanding problem in analytic number theory is to understand the size of an $L$-function at its central point. For an $L$-function $L(s)$ from the Selberg class with analytic conductor $C$ and functional equation relating values at $s$ and $1-s$, the Lindel\"{o}f hypothesis $L(\half)\ll_\epsilon C^{\epsilon}$ is expected for any $\epsilon>0$. Given an average version of the Ramanujan conjecture (which in many cases is available by the works of Iwaniec \cite{iwa}, Molteni \cite{mol}, and Xiannan Li \cite{xli}), it only requires the functional equation to prove the so called convexity bound $L(\half)\ll_{\epsilon, d} C^{\frac{1}{4}+\epsilon}$, where $d$ denotes the degree of $L(s)$. This (or a refinement of this due to Heath-Brown \cite{h-b}) is considered the trivial bound and was the best known in general until Soundararajan \cite{sou} recently proved, assuming the Ramanujan conjecture, that $L(\half)\ll_{\epsilon, d} C^{\frac{1}{4}}(\log C)^{-1+ \epsilon}$. The subconvexity problem is to save a power of $C$; that is to prove that $L(\half)\ll_{\epsilon,d} C^{\frac{1}{4}-\delta}$ for some $\delta>0$. For $L$-functions of degree 1 or 2, the problem is completely solved. This involves the work of many authors, but the contribution of Friedlander and Iwaniec is particularly noteworthy. They invented the amplifier method \cite{friiwa, iwa}, which has been used to solve many cases of  the subconvexity problem. For higher degree $L$-functions, a subconvexity bound is known only in a limited number of cases and remains a challenging and important goal. 

In this paper we study certain degree 6 $L$-functions, the Rankin-Selberg $GL(3)\times GL(2)$ $L$-functions. In a recent breakthrough, Xiaoqing Li \cite{li1} proved a subconvexity bound for the $L$-function of a self-dual Hecke-Maass form for $SL_3(\mathbb{Z})$ twisted by a Hecke-Maass form for $SL_2(\mathbb{Z})$,  or by a holomorphic Hecke cusp form for $SL_2(\mathbb{Z})$, in the eigenvalue aspect, or respectively in the weight aspect, of the $GL(2)$ form. As a corollary she derived subconvexity for a self dual $GL(3)$ form on the critical line. Blomer \cite{blom} considered this problem in the level aspect and proved subconvexity for $GL(3)\times GL(2)$ $L$-functions where the twist is by special Hecke-Maass forms of prime square level. For prime level however, subconvexity is still unknown and this is the problem which we consider.

Let $H_k^{\star}(q)$ denote the set of holomorphic cusp forms of weight $k$ which are newforms of level $q$ with trivial nebentypus in the sense of Atkin-Lehner Theory \cite{atk}. Fix $g$ a self-dual Hecke-Maass form for $SL_3(\mathbb{Z})$ which is unramified at infinity.  Let  $L(s, g\times f)$ denote the Rankin-Selberg convolution of $g$ with $f\in H_k^{\star}(q)$. Kim and Shahidi \cite{kimshah} have shown that this is in fact an automorphic $L$-function.  We normalize to have the central point at $s=\half$. The analytic conductor in the $q$ aspect equals $q^3$, so that the convexity bound is $q^{\frac{3}{4}+\epsilon}$.

In the works of Xiaoqing Li and Blomer, a study of the first moment of the $L$-function at $s=\half$ is enough to yield subconvexity. For example, in the weight aspect, the analytic conductor of $L(s, g\times f)$ equals $k^6$ so that the convexity bound is $k^{\frac{3}{2}+\epsilon}$. We further know by a result of Lapid \cite{lap} that $L(g\times f, \half)\ge 0$. Hence if we had the expected (by the Lindel\"{o}f hypothesis) upper bound
\begin{align}
\label{kaspect} \sum_{f\in H_k^{\star}(q)} L(\half, g\times f)\ll_{q,\epsilon} k^{1+\epsilon},
\end{align}
then dropping all but one term of this sum would immediately yield subconvexity. Although Li does not establish (\ref{kaspect}), she studies a similar first moment with an extra averaging over $k$. On the other hand, if we had the expected upper bound
\begin{align}
\sum_{f\in H_k^{\star}(q)}  L(\half, g\times f)\ll_{k,\epsilon} q^{1+\epsilon},
\end{align}
then dropping all but one term would not yield any useful bound in the $q$ aspect. One could try to estimate the second moment, but this seems difficult. Thus a new ingredient is needed. We make use of an amplifier to prove
\begin{theorem}\label{main}
Fix $k> 10^{6}$ an even number. Let $q$ be a prime number. Suppose that for $f_0\in H^{\star}_k(q)$ we have
\begin{align}
\label{mainassumption} \sum_{n<L} a_{f_0}(n)^2\gg_\epsilon L^{1-\epsilon}
\end{align}
for $L>q^{1/4+1/2001}$, where $a_{f_0}(n)$ is the $n$-th Fourier coefficient of $f_0$ as defined in (\ref{fouexp}).
Then
\begin{align}
\label{mainclaim} L(\half, g\times f_0)\ll q^{3/4-1/2001}.
\end{align}
\end{theorem}
\noindent One of the fundamental contributions to the subconvexity problem for degree 2 $L$-functions is Iwaniec's conditional proof of subconvexity for Hecke-Maass $L$-functions in the eigenvalue aspect \cite{iwa}, in which an assumption just like (\ref{mainassumption}) is made. The assumption in the theorem is expected of all $f_0 \in H^{\star}_k(q)$ and any $L>q^{\epsilon}$. It is known to be true for almost all $f_0 \in H^{\star}_k(q)$ and would follow, for instance, from a strong subconvexity bound for $L(\frac{1}{2}+it, f_0\times f_0)$ in the $q$ aspect. It is interesting that a bound on one $L$-function can imply a bound on another very different one. 

The exponent in our subconvexity bound and the lower bound for $k$ are not optimal. We have concentrated on a method to break the convexity bound, leaving the task of finding the best parameters to a time when the theorem can be made unconditional. 

\subsection{$L$-functions}

Every $f\in H_k^{\star}(q)$ has a Fourier expansion of the type
\begin{align}
\label{fouexp} f(z) = \sum_{n=1}^{\infty} a_f(n) n^{\frac{k-1}{2}} e(nz)
\end{align}
for $\Im z>0$, where $e(z)=e^{2\pi i z}$, $a_f(n)\in \mathbb{R}$ and $a_f(1)=1$. The coefficients $a_f(n)$ satisfy the multiplicative relation
\begin{align}
\label{hmult} a_f(n)a_f(m) = \sum_{\substack{d|(n,m)\\(d,q)=1}} a_f\Big(\frac{nm}{d^2}\Big)
\end{align}
and Deligne's bound $a_f(n)\le d(n)\ll n^{\epsilon}$. Here and throughout the paper, $\epsilon$ denotes an arbitrarily small positive constant, but not necessarily the same one from one occurrence to the next, and any implied constant may depend on $\epsilon$. Also, $q$ will always be a prime number. The $L$-function associated to $f$ is defined as
\begin{align}
L(s,f) = \sum_{n=1}^{\infty} \frac{a_f(n)}{n^s}
\end{align}
for $\Re(s)>1$. This satisfies the functional equation
\begin{multline}
\label{funct1} q^{\frac{s}{2}}    \pi^{-s}  \Gamma\Big( \frac{s+ \frac{k-1}{2}}{2} \Big)\Gamma\Big( \frac{s+ \frac{k+1}{2}}{2} \Big) L(s,f) \\
= \epsilon_f q^{\frac{1-s}{2}}  \pi^{-s}  \Gamma\Big( \frac{1-s+ \frac{k-1}{2}}{2} \Big)\Gamma\Big( \frac{1-s+ \frac{k+1}{2}}{2} \Big)  L(1-s,f),
\end{multline}
where 
\begin{align}
\label{epsf} \epsilon_f = - i^k  a_f(q) q^{\frac{1}{2}}= \pm 1.
\end{align}
The left hand side of (\ref{funct1}) analytically continues to an entire function. The facts above can be found in \cite{iwaniec}.

We fix a self-dual Hecke-Maass form $g$ of type $(\nu,\nu)$ for $SL_3(\mathbb{Z})$. We refer to \cite{gold}, especially Chapter 6, and follow its notation. We write $A(n,m)=A(m,n)$ for the Fourier coefficients of $g$ in the Fourier expansion (6.2.1) of \cite{gold}, normalized so that $A(1,1)=1$. The $L$-function associated to $g$ is defined as
\begin{align}
\label{l(g)} L(s,g) = \sum_{n=1}^{\infty} \frac{A(n,1)}{n^s}
\end{align}
for $\Re(s)>1$. The coefficients $A(n,1)$ are real. $L(s,g)$ is actually the symmetric-square $L$-function of a Hecke-Maass form for $SL_2(\mathbb{Z})$, by the work of Soudry \cite{soud}.
This implies, by the work of Kim and Sarnak \cite{kimsar}, that
\begin{align}
\label{rama} |A(n,1)| \ll n^{7/32+\epsilon}
\end{align}
and, by the work of Selberg, that
\begin{align}
\label{g1} \Re(3\nu-1)= 0.
\end{align}

We  have the Hecke relation
\begin{align}
\label{hecke3} A(n,m) = \sum_{d|(n,m)} \mu(d) A\Big(\frac{n}{d},1\Big)A\Big(1,\frac{m}{d}\Big),
\end{align}
and if $(n_1m_1,n_2m_2)=1$, we have
\begin{align}
A(n_1n_2, m_1m_2)= A(n_1,m_1)A(n_2,m_2).
\end{align}
By (\ref{rama}) and (\ref{hecke3}) we have
\begin{align}
A(n,m)\ll (nm)^{7/32+\epsilon}.
\end{align}
By (\ref{rama}) and Rankin-Selberg theory we have (cf. \cite{blom} for a proof):
\begin{align}
\sum_{n\le x} |A(na,b)|^2 \ll x(ab)^{7/16+\epsilon}.
\end{align}
This together with (\ref{hecke3}) and the Cauchy-Schwarz inequality yields
\begin{align}
\label{avg} \sum_{\substack{n<x\\m<y}} |A(na,mb)| \ll (xy)^{1+\epsilon} (ab)^{7/32 +\epsilon}.
\end{align}

The Rankin-Selberg $L$-function $L(s, g\times f)$ is defined as
\begin{align}
L(s, g\times f)= \sum_{n,r\ge 1} \frac{A(r,n)a_f(n)}{(r^2n)^s}
\end{align}
for $\Re(s)>1$. It satisfies the functional equation
\begin{align}
\label{funct2}  q^{\frac{3s}{2}} G(s)  L(s, g\times f) = \epsilon_{g\times f} q^{\frac{3(1-s)}{2}}  {G}(1-s)  L(1-s, g\times {f}),
\end{align}
where $\epsilon_{g\times f} = (\epsilon_f)^3=\epsilon_f$ and
\begin{multline}
G(s) = \pi^{-3s} \Gamma\Big( \frac{s+\frac{k+1}{2} + 3\nu-1}{2} \Big)\Gamma\Big( \frac{s+\frac{k+1}{2} }{2} \Big)\Gamma\Big( \frac{s+\frac{k+1}{2} +1-3\nu}{2} \Big)\\
 \times \Gamma\Big( \frac{s+ \frac{k-1}{2} + 3\nu-1}{2} \Big)\Gamma\Big( \frac{s+ \frac{k- 1}{2}}{2} \Big)\Gamma\Big( \frac{s+\frac{k-1}{2} +1-3\nu}{2} \Big).
\end{multline}
The left hand side of (\ref{funct2}) analytically continues to an entire function. To study $L(1/2,g\times f)$, we first express it as a weighted Dirichlet series.

\begin{lemma}{\bf Approximate functional equation} \label{afelemma}

Let 
\begin{align}
\label{v1def} V(x)= \frac{1}{2\pi i} \int_{(\sigma)} x^{-s} \frac{G(\half + s)}{G(\half)} \frac{ds}{s}\end{align}
for $x, \sigma>0$.
We have
\begin{multline}
\label{afe1} L(\half,f\times g) = \sum_{n,r\ge 1} \frac{a_f(n)A(r,n)}{r\sqrt{n}}V\Big(\frac{r^2n}{q^{3/2+1/300}}\Big)\\+ \epsilon_{g\times f} \sum_{n,r\ge 1} \frac{a_f(n)A(r,n)}{r\sqrt{n}}{V}\Big(\frac{r^2n}{q^{3/2-1/300}}\Big).
\end{multline}

For any $A>0$ and integer $B\ge 0$ we have that
\begin{align}
\label{vbound} V^{(B)}(x) \ll_B x^B (1+x)^{-A},
\end{align}
so that the first sum in (\ref{afe1}) is essentially supported on  $r^2n<q^{3/2+1/300+\epsilon}$ and the second sum is essentially supported on  $r^2n<q^{3/2-1/300+\epsilon}$.
\end{lemma}
\proof
The proof of this standard result may be found in Theorem 5.3 of \cite{iwakow}. 
\endproof
\noindent Note that the sums in (\ref{afe1}) are of different lengths. This will result in less work with the second sum, which contains the root number $\epsilon_{g\times f}$. 

\subsection{Trace formula}

We have Weil's estimate for the Kloosterman sum:
\begin{align}
\label{weil} |S(n,m;c)| = \Big| \summ_{h\bmod c} e\Big(\frac{nh+m\overline{h}}{c}\Big)\Big| \le (n,m,c)^{\frac{1}{2}} c^{\frac{1}{2}} d(c).
\end{align}
Here $\star$ denotes that the summation is restricted to $(h,c)=1$, and $\overline{h}h\equiv 1 \bmod c$.
We will also need the following estimates for the $J$-Bessel function (see \cite{gr} and \cite{wat}). For $x>0$ we have
\begin{align}
\label{jbound} J_{k-1}(x)\ll \min(x^{k-1},x^{-1/2}).
\end{align}
For $x>0$ and integers $B>0$ we have
\begin{align}
\label{jbound2} J_{k-1}^{(B)}(x)\ll_{B} x^{-B}+x^{-1/2}.
\end{align}

For any complex numbers $\alpha_f$, define the weighted sum
\begin{align}
\sumh_{f\in H_k^{\star}(q)}  \alpha_f = \sum_{f\in H_k^{\star}(q)} \frac{\alpha_f}{\zeta(2)^{-1} L(1, \text{sym}^2 f)},
\end{align}
where $ L(s, \text{sym}^2 f)$ denotes the symmetric-square $L$-function of $f$. The arithmetic weights above occur naturally in the Petersson trace formula (\ref{ptf}) and the following trace formula for newforms. Define
\begin{align}
\Delta^{\star}_{k,q}(n,m)=\frac{12}{q(k-1)}& \sumh_{f\in H_k^{\star}(q)} a_f(n)a_f(m).
\end{align}

\begin{lemma}{\bf Trace formula.} \label{tf}

(i) We have
\begin{align}
\label{ptf} \Delta^{\star}_{k,1}(n,m)= \delta(n,m)  + 2\pi i^k \sum_{c\ge 1} \frac{S(n,m;c)}{c}J_{k-1}\Big(\frac{4\pi\sqrt{nm}}{c}\Big),
\end{align}
where $\delta(n,m)$ equals $1$ if $n=m$ and $0$ otherwise.

(ii) Let $q$ be a prime. If $(m,q)=1$ and $q^2\nmid n$ then
\begin{multline}
 \label{trace}\Delta^{\star}_{k,q}(n,m)  = \delta(n,m) + 2\pi i^k \sum_{c\ge 1} \frac{S(n,m;cq)}{cq}J_{k-1}\Big(\frac{4\pi\sqrt{nm}}{cq}\Big)\\
 -\frac{1}{q[\Gamma_0(1):\Gamma_0((n,q))]}\sum_{i=0}^{\infty} q^{-i} \Delta^{\star}_{k,1}(n, mq^{2i}).
 \end{multline}
\end{lemma}
\proof
(\ref{ptf}) can be found in \cite{iwaniec}. See Proposition 2.8 of \cite{ils} for (\ref{trace}).
\endproof
\noindent Note that if $q|n$ then the last line of (\ref{trace}) is $\ll q^{-2+\epsilon}(nm)^{\epsilon}$ since $[\Gamma_0(1):\Gamma_0(q)]>q$.

\subsection{Voronoi summation}

The GL(3) Voronoi summation formula (\ref{voron}) was found by Miller and Schmid \cite{milsch}. Goldfeld and Li \cite{golli} later gave another proof. 

\begin{lemma}{\bf GL(3) Voronoi Summation.}\label{vsf}
Let $\psi$ be a smooth, compactly supported function on the positive real numbers and $(b,d)=1$. We have
\begin{align}
\label{voron} \sum_{n\ge 1} A(r,n) e\Big(\frac{n\overline{b}}{d}\Big) \psi\left( \frac{n}{N} \right) =
 \sum_{\pm} \frac{\pi^{\frac{3}{2}}}{2} d\sum_{\substack{n\ge 1\\ l | dr}} \frac{A(n,l)}{nl}S\Big(rb,\pm n;\frac{dr}{l}\Big)\Psi^{\pm}\Big(\frac{n}{d^3N^{-1}rl^{-2}}\Big),
\end{align}
where we define
\begin{align}
&\Psi^{\pm}(X)=  X \frac{1}{2\pi i}\int_{(\sigma)} (\pi^3 X)^{-s} H^{\pm}(s) \widetilde{\psi}(1-s) ds ,\\
&\nonumber H^{\pm}(s)= \frac{\Gamma\Big(\frac{s+3\nu-1}{2}\Big)\Gamma\Big(\frac{s}{2}\Big)\Gamma\Big(\frac{s+1-3\nu}{2}\Big)}{\Gamma\Big(\frac{1-s-1+3\nu}{2}\Big)\Gamma\Big(\frac{1-s}{2}\Big)\Gamma\Big(\frac{1-s+3\nu-1}{2}\Big)} 
\mp i \frac{\Gamma\Big(\frac{1+s+3\nu-1}{2}\Big)\Gamma\Big(\frac{1+s}{2}\Big)\Gamma\Big(\frac{1+s+1-3\nu}{2}\Big)}{\Gamma\Big(\frac{2-s+1-3\nu}{2}\Big)\Gamma\Big(\frac{2-s}{2}\Big)\Gamma\Big(\frac{2-s+3\nu-1}{2}\Big)}
\end{align}
for $\sigma > 0$, where $ \widetilde{\psi}$ denotes the Mellin transform of $\psi$.
\end{lemma}
\noindent Writing $s=\sigma+it$, by Stirling's approximation of the gamma function we have
\begin{align}
\label{Hbound} H^{\pm}(s)\ll_\sigma (1+|t|)^{3\sigma}.
\end{align}

\subsection{Amplifier method}

Let $f_0\in H_k^{\star}(q)$. Define the amplifier
\begin{align}
A(f)= \sum_{n< q^{1/4+1/2000}} \frac{a_{f_0}(n)a_f(n)}{\sqrt{n}}.
\end{align}
The assumption (\ref{mainassumption}) implies that $A(f)$ is `amplified' at $f=f_0$:
\begin{align}
\label{hyp} \sum_{n<q^{1/4+1/2000}} \frac{a_{f_0}(n)^2}{\sqrt{n}} \gg q^{1/8+1/4000-\epsilon}.
\end{align}
This can be seen by partial summation together with the following upper bound given in \cite{mol}: 
\begin{align}
\sum_{n<L} a_{f_0}(n)^2 \ll L^{1+\epsilon}
\end{align}
for all $L>q^{\epsilon}$.
Theorem \ref{main} will be deduced from the following.
\begin{proposition}\label{prop1}
We have
\begin{align}
\sumh_{f\in H_k^{\star}(q)} L(\half, g\times f) A(f)^2 \ll q^{1+\epsilon}.
\end{align}
\end{proposition}
\noindent By Lapid's work, we have that $ L(\half, g\times f)\ge 0$. Now if we drop all but the term corresponding to $f_0$ then we have
\begin{align}
\frac{1}{L(1, \text{sym}^2 f_0)}L(\half, g\times f_0) A(f_0)^2\ll q^{1+\epsilon}.
\end{align}
Using the trivial bound $L(1, \text{sym}^2 f_0)\ll q^{\epsilon}$ and the assumption (\ref{hyp}), the subconvexity bound (\ref{mainclaim}) follows.

By (\ref{hmult}) we may write
\begin{align}
A(f)^2=\sum_{m< q^{1/2+1/1000}} \frac{x_m a_f(m)}{\sqrt{m}},
\end{align}
for some numbers $x_m\ll q^{\epsilon}$. By (\ref{avg}), Proposition \ref{prop1} follows from
\begin{proposition}\label{prop2}
Let $m<q^{1/2+1/1000}$ be a natural number. We have
\begin{align}
\label{p2} \frac{1}{q} \sumh_{f\in H_k^{\star}(q)} L(\half, g\times f) a_f(m) \ll \frac{q^{\epsilon}}{\sqrt{m}}\Big(1+\sum_{r< q^2} \frac{|A(r,m)|}{r}\Big).
\end{align}
\end{proposition}

\section{Proof of Proposition \ref{prop2}}

In this section we reduce the proof of Proposition \ref{prop2} to two claims. By Lemma \ref{afelemma} and (\ref{hmult}), we have that the left hand side of (\ref{p2}) is
\begin{multline}
\label{t1} \ll \sum_{n,r\ge 1} \frac{A(r,n)}{r\sqrt{n}}V\Big(\frac{r^2n}{q^{3/2+1/300}}\Big)\Delta^{\star}_{k,q}(n,m)\\
+q^{1/2} \sum_{n,r\ge 1} \frac{A(r,n)}{r\sqrt{n}}{V}\Big(\frac{r^2n}{q^{3/2-1/300}}\Big)\Delta^{\star}_{k,q}(nq,m).
\end{multline}
By Lemma \ref{tf}, the first line of (\ref{t1}) is
\begin{multline}
\label{t2} \ll \sum_{r<q^2} \frac{|A(r,m)|}{r\sqrt{m}} + \Big|\sum_{n,r\ge 1} \frac{A(r,n)}{r\sqrt{n}}V\Big(\frac{r^2n}{q^{3/2+1/300}}\Big)\sum_{c\ge 1} \frac{S(n,m;cq)}{cq} J_{k-1}\Big(\frac{4\pi\sqrt{nm}}{cq}\Big)\Big|
\\+q^{-1+\epsilon} \sum_{f'\in H_k^{\star}(1)}\Big| \sum_{n,r\ge 1} \frac{A(r,n)a_{f'}(n)}{r\sqrt{n}}V\Big(\frac{r^2n}{q^{3/2+1/300}}\Big)\Big|.
\end{multline}
Using (\ref{v1def}), the last sum over $n$ and $r$ above can be written as an integral involving $L(s,g\times f')$. The line of integration can be moved to $-\infty$ to see that the last line of (\ref{t2}) is $\ll q^{-1+\epsilon}$. We will prove
\begin{lemma}\label{mainlem1}
Let $m<q^{1/2+1/1000}$. We have
\begin{align}
\label{lem1} \sum_{n,r\ge 1} \frac{A(r,n)}{r\sqrt{n}}V\Big(\frac{r^2n}{q^{3/2+1/300}}\Big)\sum_{c\ge 1} \frac{S(n,m;cq)}{cq} J_{k-1}\Big(\frac{4\pi\sqrt{nm}}{cq}\Big) \ll \frac{q^{\epsilon}}{\sqrt{m}}.
\end{align}
\end{lemma}

For the second line of (\ref{t1}), we first consider the contribution of the terms with $(n,q)=1$. By Lemma \ref{tf} and the remark immediately following, the contribution of such terms is
\begin{align}
\ll q^{1/2} \sum_{\substack{n,r\ge 1\\(n,q)=1}} \frac{A(r,n)}{r\sqrt{n}}{V}\Big(\frac{r^2n}{q^{3/2-1/300}}\Big)\sum_{c\ge 1} \frac{S(nq,m;cq)}{cq} J_{k-1}\Big(\frac{4\pi\sqrt{nqm}}{cq}\Big)+O(q^{-1/2}).
\end{align}
In the sum above, the contribution of the terms with $c>q^{1/2}$ is $\ll q^{-100}$. Thus we may assume that $(c,q)=1$, so that we have $S(nq,m,cq)= -S(n\overline{q},m,c)$. We may extend the sum to all natural numbers $n$, with an error of
\begin{align}
\ll \sum_{n,r\ge 1} \frac{A(r,nq)}{r\sqrt{n}}{V}\Big(\frac{r^2n}{q^{1/2-1/300}}\Big)
\frac{1}{q}\sum_{c\ge 1} \frac{S(n,m;c)}{c} J_{k-1}\Big(\frac{4\pi\sqrt{nm}}{c}\Big)\ll q^{-1/2},
\end{align}
on observing that the $c$-sum equals $\Delta^{\star}_{k,1}(n,m)\ll q^{\epsilon}$ and using (\ref{avg}). We will prove
\begin{lemma}\label{mainlem2}
Let $m<q^{1/2+1/1000}$. We have
\begin{align}
q^{1/2} \sum_{n,r\ge 1} \frac{A(r,n)}{r\sqrt{n}}{V}\Big(\frac{r^2n}{q^{3/2-1/300}}\Big)\sum_{c\ge 1} \frac{S(n\overline{q},m;c)}{cq} J_{k-1}\Big(\frac{4\pi\sqrt{nmq}}{cq}\Big)\ll q^{-1}.
\end{align}
\end{lemma}
\noindent Finally we must consider the terms of the second line of (\ref{t1}) with $q|n$. The contribution of these terms is
\begin{align}
\ll \sum_{n,r\ge 1} \frac{A(r,nq)}{r\sqrt{n}}V\Big(\frac{r^2n}{q^{1/2-1/300}}\Big)\Delta^{\star}_{k,q}(nq^2,m) \ll q^{-1/2},
\end{align}
on using (\ref{hmult}) and (\ref{epsf}) to see that $\Delta^{\star}_{k,q}(nq^2,m)=q^{-1} \Delta^{\star}_{k,q}(n,m)$ and then using (\ref{avg}) to bound the sum absolutely. This completes the proof of Proposition \ref{prop2}.

\subsection{Sketch.}

Before starting on the proofs of the lemmas presented in this section, we give a rough sketch of the argument for the main lemma, Lemma \ref{mainlem1}. Since $m\approx q^{1/2}$ and $r^2n$ is essentially bounded by about $q^{3/2}$ in (\ref{lem1}), the value of the $J$-Bessel function will be very small unless $n\approx q^{3/2}$, $r\approx q^{\epsilon}$ and $c\approx q^{\epsilon}$. Consider the range $q^{3/2}<n<2q^{3/2}$, $r=1$ and $c=1$. Opening the Kloosterman sum, a part of what we must bound in (\ref{lem1}) is
\begin{align}
\summ_{h \bmod q} e(mh/q) \sum_{q^{3/2}< n< 2q^{3/2}} A(1,n) e(n\overline{h}/q).
\end{align}
The weight function $V$ and the $J$-Bessel function have been ignored because they are roughly constant in this range. We apply the $GL(3)$ Voronoi summation formula to exchange the $n$-sum for another sum of length about $q^3/q^{3/2}=q^{3/2}$. A part of what we must bound is then
\begin{align}
\summ_{h \bmod q} e(mh/q) \sum_{q^{3/2}< n< 2q^{3/2}} A(n,1) S(n,h,q)
\end{align}
We have $\summ_{h \bmod q} e(mh/q) S(n,h,q) \approx q e(n\overline{m}/q)$. By reciprocity (the Chinese Remainder Theorem), we have $e(n\overline{m}/q)=e(n/mq)e(-n\overline{q}/m)\approx e(-n\overline{q}/m)$, since $n\approx q^{3/2} \approx mq$. Thus we must bound
\begin{align}
 \sum_{q^{3/2}< n< 2q^{3/2}} A(n,1) e(-n\overline{q}/m).
\end{align}
The new modulus $m$ of the exponential is much smaller than the original modulus $q$. We apply the $GL(3)$ Voronoi summation formula once again, to exchange the $n$-sum for another sum of length about $m^3/q^{3/2}\approx 1$. We must bound
\begin{align}
\sum_{ n< q^{\epsilon}} A(1,n) S(-n,q,m).
\end{align}
We bound this sum absolutely, using Weil's bound for the Kloosterman sum.

\section{Proof of Lemma \ref{mainlem2}}

Write $S(n\overline{q},m,c) = \summ_{h \bmod c} e(n\overline{qh}/c)e(mh,c)$. As noted above, by (\ref{jbound}) we may assume that $c<q^{1/2}$. By (\ref{vbound}), we may also assume that $r^2<q^{3/2}$. Thus to prove Lemma \ref{mainlem2}, it is enough to show that
\begin{align}
\sum_{n} A(r,n)e(n\overline{qh}/c)n^{-1/2}J_{k-1}\Big(\frac{4\pi\sqrt{nmq}}{cq}\Big) {V}\Big(\frac{r^2n}{q^{3/2-1/300}}\Big) \ll q^{-2}.
\end{align}
We consider this sum in dyadic intervals. For $N>0$, let $\omega(x)$ be a smooth function, compactly supported on $[1,2]$ and satisfying $\omega^{(B)}(x)\ll_B 1$ and let
\begin{align}
W(x)= x^{-1/2}  J_{k-1}\Big( \frac{4\pi\sqrt{xNmq^{-1}}}{c} \Big)  V\Big(\frac{xr^2N}{q^{3/2-1/300}}\Big)\omega(x).
\end{align}
It is enough to show that
\begin{align}
\label{esum} \sum_{n}A(r,n)e(n\overline{qh}/c)W\Big(\frac{n}{N}\Big) \ll q^{-3}
\end{align}
for
\begin{align}
\label{r2n} r^2N<q^{3/2-1/300+\epsilon}
\end{align}
and 
\begin{align}
\label{jc} \frac{\sqrt{Nmq^{-1}}}{c} > q^{-1/10^5}.
\end{align}
We enforce the conditions (\ref{r2n}) and (\ref{jc}) since otherwise (\ref{esum}) follows easily by (\ref{vbound}) and (\ref{jbound}).

Applying Lemma {\ref{vsf} to (\ref{esum}), it is enough to show that
\begin{align}
\label{fsum} \sum_{\substack{n\ge 1\\ l|rc}} \frac{A(n,l)}{nl} S(rhq,\pm n,rc/l)  \mathcal{W}\Big(\frac{nNl^2}{c^3r}\Big)  \ll q^{-4},
\end{align}
where
\begin{align}
\label{rpsi} \mathcal{W} (X)= \int_{(\sigma)} X^{1-s} H^{\pm}(s)\widetilde{W}(1-s)  ds,
\end{align}
for $\sigma>0$ and
\begin{align}
\widetilde{W}(1-s)= \int_1^2 x^{-s} W(x) dx.
\end{align}
By (\ref{jbound2}) and (\ref{jc}) we have for $B>0$,
\begin{align}
 W^{(B)}(x) \ll_B  \Big( \frac{\sqrt{Nmq^{-1}}}{c} \Big)^{B}q^{B/10^{5}}.
\end{align}
Thus, writing $s=\sigma+it$, we have by integration by parts $B$ times,
\begin{align}
\label{bou} \widetilde{W}(1-s) \ll_B |t|^{-B}  \Big( \frac{\sqrt{Nmq^{-1}}}{c} \Big)^{B}q^{B/10^5}.
\end{align}
Using this bound with $B=\lfloor 3\sigma +5 \rfloor$ and (\ref{Hbound}), we have
\begin{align}
\mathcal{W}\Big(\frac{nNl^2}{c^3r}\Big) &\ll_\sigma q^{10} \Big( \frac{nl^2N}{c^3r} \Big)^{-\sigma}  \Big( \frac{\sqrt{Nmq^{-1}}}{c} \Big)^{3\sigma} q^{3\sigma/10^5}\\
\nonumber &\ll q^{10} (nl^2)^{-\sigma} \Big( \frac{r^2Nm^3}{q^{3}} \Big)^{\sigma/2}q^{3\sigma/10^5}. 
\end{align}
By (\ref{r2n}) and the assumptions of the lemma, we have $r^2Nm^3 \ll q^{3-1/10^4}$. Thus taking $\sigma$ large enough proves (\ref{fsum}).

\section{Proof of Lemma \ref{mainlem1}}

We consider the $n$-sum in dyadic intervals. For $N_1>0$, let $\omega_1(x)$ be a smooth function, compactly supported on $[1,2]$ and satisfying $\omega_1^{(B)}(x)\ll_B 1$ and let
\begin{align}
W_1(x)= x^{-1/2}  J_{k-1}\Big( \frac{4\pi\sqrt{xN_1mq^{-2}}}{c} \Big)  V\Big(\frac{xr^2N_1}{q^{3/2+1/300}}\Big)\omega_1(x).
\end{align}
It is enough to show that
\begin{align}
\label{ls0} \sum_{n\ge 1} A(r,n)S(n,m;cq)W_1\Big(\frac{n}{N_1}\Big)\ll \frac{q^{1+\epsilon}\sqrt{N_1}}{\sqrt{m}}
\end{align}
for
\begin{align}
q^{3/2-1/990}<&N_1<q^{3/2+1/300+\epsilon},\\
\nonumber q^{1/2-1/290} <&m < q^{1/2+1/1000},\\
\nonumber &r< q^{1/450},\\
\nonumber &c< q^{1/450}.
\end{align}
We may assume the conditions above, since otherwise (\ref{ls0}) follows easily by (\ref{vbound}) and (\ref{jbound}). Thus it is enough to prove that
\begin{align}
\label{ls1} \sum_{n\ge 1} A(r,n)S(n,m;cq)W_1\Big(\frac{n}{N_1}\Big)\ll q^{3/2-1/990}.
\end{align}

We apply Lemma \ref{vsf} to the left hand side of (\ref{ls1}) after writing $S(n,m;cq)=\summ_{h \bmod cq} e((n\overline{h}+mh)/cq)$. We need to show that
\begin{align}
\label{off3}  cq \summ_{h \bmod cq} e(mh/cq) \sum_{\substack{n\ge1\\ l | cqr}} \frac{ A(n, l)}{nl} S(rh, n; qcr/l) \mathcal{W}_1\Big(\frac{nN_1l^2}{c^3q^3r}\Big)\ll  q^{3/2-1/990},
\end{align}
where
\begin{align}
\label{defw2} \mathcal{W}_1 (X)= X \int_{(\sigma)} (\pi^3X)^{-s} H^\pm(s) \widetilde{W_1}(1-s) ds
\end{align}
for $\sigma>0$.
Note that $W_1^{(B)}(x)\ll_B q^{B/450}$ so that by integrating by parts $B$ times we have for $s=\sigma+it$,
\begin{align}
\label{w1bounds} \widetilde{W_1}(1-s) = \int_{1}^{2} x^{-s} W_1(x) dx \ll_B |t|^{-B} q^{B/450}.
\end{align}
We can use this bound together with (\ref{Hbound}) to estimate $\mathcal{W}_1(X)$. If $X>q^{1/150+\epsilon}$, we can take $\sigma$ in (\ref{defw2}) to be very large and $B=\lfloor 3\sigma +5 \rfloor$ in (\ref{w1bounds}) to see that $\mathcal{W}_1(X)\ll X^{-2}q^{-100}$. If $ X \le q^{1/150+\epsilon}$, we take $\sigma=\epsilon$  in (\ref{defw2}) and $B=2$ in (\ref{w1bounds}) to see that 
\begin{align}
\label{w2bound} \mathcal{W}_1(X)\ll q^{1/200}X.
\end{align}
So the $n$-sum in (\ref{off3}) is essentially supported on $n<\frac{q^{3+1/150+\epsilon}c^{3}r}{N_1 l^2}<q^{3/2+1/60}$.

The contribution to (\ref{off3}) by the terms with $q|l$ is negligible since if $q|l$ then $\frac{nN_1l^2}{c^3q^3r}>q^{1/150+\epsilon}$ and we have just seen that then $\mathcal{W}_1\Big(\frac{nN_1l^2}{c^3q^3r}\Big)\ll q^{-100}$.
Henceforth fix $l|cr$, so that $l< q^{1/225}$.

We open the Kloosterman sum: $S(rh, n; qcr/l) = \summ_{u \bmod qcr/l} e((rhu+n\overline{u})l/qcr)$. For (\ref{off3}), it is enough to show that
\begin{multline}
 \label{off4} \summ_{u \bmod qcr/l}  \sum_{n\ge 1} \frac{A(n,l) e(n\overline{u}l/qcr)}{n}\mathcal{W}_1\Big(\frac{nNl^2}{c^3q^3r}\Big) \summ_{h \bmod cq} e(h(ul+m)/cq)\ll q^{1/2-1/300}.
\end{multline}
The innermost sum above, a Ramanujan sum, equals
\begin{align}
\Big(\summ_{h \bmod q} e(h(ul+m)/q)\Big)\Big(\summ_{h \bmod c} e(h(ul+m)/c)\Big),
\end{align}
since $(c,q)=1$. 
Note that $\summ_{h \bmod q} e(h(ul+m)/q)$ equals $-1$ or $q-1$ according as $ul \not\equiv m\bmod q$ or $ul \equiv m \bmod q$ respectively. So the left hand side of (\ref{off4}) equals
\begin{align}
& - \sum_{n\ge 1}  \frac{A(n,l)}{n}\mathcal{W}_1\Big(\frac{nNl^2}{c^3q^3r}\Big)  \summ_{h \bmod c}  e(hm/c) S(n,qhr;qcr/l) \label{off5} \\
& +\ q\sum_{n\ge 1} \frac{A(n,l) }{n}\mathcal{W}_1\Big(\frac{nNl^2}{c^3q^3r}\Big) \summ_{h \bmod c}  e(hm/c) \summ_{\substack{u \bmod qcr/l\\ u\equiv m \overline{l}\bmod q}}e\Big(\frac{uqhr+\overline{u}n}{qcr/l} \Big).\nonumber
\end{align}
Since $(cr/l,q)=1$, we have $S(n,qhr;qcr/l)=S(n\overline{q},hr;cr/l)S(0,n;q)$. This product of a Kloosterman sum and a Ramanujan sum is $\ll q^{1+1/225}$ if $q|n$ and $\ll q^{1/225}$ otherwise. In any case, using (\ref{avg}), the first line of (\ref{off5}) satisfies the bound required in (\ref{off4}). Now consider the second line. By the Chinese Remainder Theorem, for $(u,qcr/l)=1$ and $u\equiv m\overline{l} \bmod q$, we can write $u=m\overline{cr}(cr/l) + vq$, where $cr\overline{cr}\equiv 1 \bmod q$ and $(v,cr/l)=1$.  We have 
\begin{align}
&e\Big(\frac{uqhr}{qcr/l}\Big)=e\Big(\frac{vhrq}{cr/l}\Big),\\
&e\Big(\frac{n\overline{u}}{qcr/l}\Big)=e\Big(\frac{n\overline{ucr/l}}{q}\Big)e\Big(\frac{n\overline{uq}}{cr/l}\Big)=e\Big(\frac{nl^2\overline{mcr}}{q}\Big)e\Big(\frac{n\overline{vq^2}}{cr/l}\Big).\nonumber
\end{align}
Thus (\ref{off4}) is reduced to showing
\begin{multline}
\label{off06}  \sum_{n\ge 1}   \frac{A(n,l)}{n} e\Big(\frac{nl^2\overline{mcr}}{q}\Big)\mathcal{W}_1\Big(\frac{nNl^2}{c^3q^3r}\Big) \summ_{h \bmod c} e\Big(\frac{hm}{c}\Big)S(n\overline{q^2},hrq;cr/l)\\
\ll  q^{-1/2-1/300}.
\end{multline}

Now comes a crucial step. By the Chinese Remainder Theorem we have 
\begin{align}
e\Big(\frac{nl^2\overline{mcr}}{q}\Big)= e\Big(\frac{nl^2}{mcrq}\Big)e\Big(\frac{-nl^2\overline{q}}{mcr}\Big).
\end{align}
For $n<q^{3/2+1/60}$, the exponential factor $e\big(\frac{nl^2}{mcrq}\big)$ has amplitude at most $q^{1/30}$ and will be absorbed into the weight function.  Let $N_2<q^{3/2+1/60+\epsilon}$, let $\omega_2(x)$ be a smooth function, compactly supported on $[1,2]$ and satisfying $\omega_2^{(B)}(x)\ll_B 1$ and let
\begin{align}
W_2(x)= x^{-1} e\Big(\frac{x N_2 l^2}{mcrq}\Big) \mathcal{W}_1\Big(\frac{x N_2 N_1 l^2}{c^3q^3r}\Big)  \omega_2(x).
\end{align}
It is enough to prove that
\begin{multline}
\sum_{n\ge 1} A(n,l) e\Big(\frac{-nl^2\overline{q}}{mcr}\Big) W_2\Big( \frac{n}{N_2} \Big) \summ_{h \bmod c} e\Big(\frac{hm}{c}\Big)S(n\overline{q^2},hrq;cr/l)\\
\ll N_2 q^{-1/2-1/300}.
\end{multline}
Combining the exponential factors, it is enough to prove that
\begin{align}
\label{off7} \sum_{n\ge 1} A(n,l) e(n\overline{b}/d) W_2\Big( \frac{n}{N_2} \Big)\ll  N_2 q^{-1/2-1/100}.
\end{align}
for $d<q^{1/2+1/150}$ and (b,d)=1.
We use the Voronoi summation formula again. Applying Lemma \ref{vsf}  to the left hand side of (\ref{off7}), it suffices to show that
\begin{align}
\label{off8} d \sum_{\substack{n\ge 1\\ \ell | d l}} \frac{A(\ell,n)}{n\ell} S(lb, n ; d l/ \ell) \mathcal{W}_2 \Big( \frac{nN_2l^2}{d^3r}  \Big) \ll   N_2 q^{-1/2-1/100}.
\end{align}
where
\begin{align}
\label{w6} \mathcal{W}_2(X)= X \int_{(\sigma)} (\pi^3 X)^{-s}  H^{\pm}(s) \widetilde{W_2}(1-s) ds
\end{align}
for $\sigma> 0$.
We need to estimate $\mathcal{W}_2(X)$. To this end we first note, using (\ref{w2bound}), that
\begin{align}
W_2^{(B)}(x)\ll_B   q^{B/30} q^{1/200} \frac{N_2 N_1 l^2}{c^3q^3r} \ll q^{(B+1)/30} q^{-3/2}N_2.
\end{align}
Integrating by parts $B$ times, we have for $s=\sigma+it$ the bound
\begin{align}
\label{w7} \widetilde{W_2}(1-s) = \int_{1}^{2} x^{-s} W_2(x) dx &\ll_B |t|^{-B} q^{(B+1)/30} q^{-3/2}N_2.
 \end{align}
Now we can estimate $\mathcal{W}_2(X)$. We take $\sigma=1$ in (\ref{w6}) and $B=5$ above. By (\ref{Hbound}) and (\ref{w7}) we see that
\begin{align}
\label{wwbound}\mathcal{W}_2(X)\ll N_2 q^{-3/2+1/5}.
\end{align}
Taking $\sigma=2$, we also observe that the sum in (\ref{off8}) can be restricted to $n<q^{100}$, say, with negligible error. 

Using (\ref{wwbound}), to prove (\ref{off8}) it is enough to show that
\begin{align}
\label{last} \sum_{\substack{n< q^{100}\\ \ell | dl}} \frac{|A(\ell,n)|}{n\ell}| S(lb, n ; d l/ \ell)| \ll q^{1/2-13/60}.
\end{align}
By (\ref{weil}), we have $| S(lb, n ; d l/ \ell)| \ll (lb,n,dl/\ell)^{1/2} (dl) ^{1/2+\epsilon}\ll q^{1/4+1/100}$, since $(b,d)=1$. This establishes (\ref{last}), using (\ref{avg}).

\bibliographystyle{amsplain}

\bibliography{gl3xgl2}

\end{document}